\documentclass[11pt]{article}

\usepackage[english]{babel}
\usepackage[latin1]{inputenc}
\usepackage[T1]{fontenc}
\usepackage[usenames]{color}
\usepackage{amsmath}
\usepackage{amssymb}
\usepackage{hhline}
\usepackage{multirow}
\usepackage{colortbl}
\usepackage{enumerate}
\usepackage{array}
\usepackage{soul}
\usepackage{graphicx}
\usepackage{eurosym}
\usepackage{mathrsfs}
\usepackage{dsfont}
\usepackage{amsthm}
\usepackage{makecell}
\usepackage[titletoc]{appendix}
\usepackage[titles]{tocloft}
\usepackage{indentfirst}

\newtheorem{proposition}{Proposition}[section]

\newtheorem{definition}{Definition}[section]

\title{A tensorial-parallel Chebyshev method for a differential game theory problem}

\author{Carmelo de Castro\thanks{IMUVA, Universidad de Valladolid, Paseo de Bel\'{e}n 7, Valladolid, Spain. e-mail:carmelo.castro@uva.es},  V\'{\i}ctor Gat\'{o}n\thanks{\textit{Corresponding author}. IMUVA and Dpto. Matem\'{a}tica Aplicada, Universidad de Valladolid,  Paseo de Bel\'{e}n 7, Valladolid, Spain. e-mail:victor.gaton@uva.es} and Beatriz G\'{o}mez\thanks{Universidad de Valladolid,  Paseo de Bel\'{e}n 7, Valladolid, Spain. e-mail: beatriz.gomez.martin@estudiantes.uva.es}}

\begin{document}

\maketitle

\begin{abstract}This paper concerns the design of a multidimensional Chebyshev interpolation based method for a differential game theory problem. In continuous game theory problems, it might be difficult to find analytical solutions, so numerical methods have to be applied. As the number of players grows, this may increase computational costs due to the curse of dimensionality. To handle this, several techniques may be applied and paralellization can be employed to reduce the computational time cost. Chebyshev multidimensional interpolation allows efficient multiple evaluations simultaneously along several dimensions, so this can be employed to design a tensorial method which performs many computations at the same time. This method can also be adapted to handle parallel computation and, the combination of these techniques, greatly reduces the total computational time cost. We show how this technique can be applied in a pollution differential game. Numerical results, including error behaviour and computational time cost, comparing this technique with a spline-parallelized method are also included.

\

\textbf{Keywords:} {Transboundary pollution, Differential games, Parabolic differential equations, Chebyshev multidimensional interpolation.}
\end{abstract}

\section{Introduction}

In differential game theory (see \cite{Ba}), several agents (or players) jointly control,  through their actions, a dynamical system described by differential state equations. The actions of the agents are taken in order to maximize a particular objective function (for each player) which outcome depends on the state of the system and the actions of other players. Differential game theory is broadly employed in many areas including, for example, economics, management, engineering and operations research.

In general, it might not be easy to find explicit solutions for differential game problems, even if we restrict ourselves to a small amount of players, and numerical methods have to be employed (see \cite{Falcone} or \cite{Johnson}). If collocation methods are employed, as the number of players increases, we have to deal with the so called ``curse of dimensionality'', which might boost the computational cost of the numerical methods.

Spectral methods (see \cite{Canuto}) are a class of spatial discretizations for partial differential equations with an order of convergence that depends on the regularity of the function to be approximated. Spectral methods have been successfully employed in many fields and have been proved competitive with other alternatives, both in precision and computational time cost. For example, Chebyshev interpolation has been employed in \cite{deFrutos2} and \cite{Gab} to price financial derivatives. In \cite{Ruijter}, a Fourier cosine method is employed to solve backward stochastic differential equations. Other examples are \cite{Ortiz}, \cite{Ruijter2} or \cite{Zhangb}. In game theory and optimal control, spectral methods have also been employed. In \cite{Nikoo}, a Chebyshev pseudospectral method is employed for obtaining a numerical solution of an open-loop Nash equilibrium and in \cite{Zhang} a Spectral Galerkin method is developed.

The literature in economic and environmental problems can be divided in two categories: the papers which study the economic growth theory with spatial diffusion (for example \cite{Brito}, \cite{Camacho} or \cite{Fabbri}), and papers which deal with the spatial dimension in environmental and resource economics (for example \cite{Brock}, \cite{Camacho2} or \cite{Xepa}). Concerning transboundary pollution games specifically, \cite{Dockner} and \cite{Van} are seminal papers and a survey of the literature in that area can be found in \cite{Jor}.

The differential game that we are going to employ to test our numerical method is developed in \cite{deFrutos1}, and it corresponds to a model which combines two aspects: first, the spatial aspect to the transboundary pollution dynamic games and second, the  strategic aspects to the spatial economics, in particular to the pollution control in a spatial setting.

The paper is organized as follows. In Section \ref{pollutiong} we make a brief description of the differential pollution game, which can be found in \cite{deFrutos1}, and we present the Chebyshev interpolation based algorithm that can be employed to numerically solve the game. In Section \ref{cbmethod}, we describe several algorithms that allow an efficient valuation of the polynomials involved and we show how the method can be extended to handle parallelization. Section \ref{numres} gives some numerical results, including both numerical error behaviour and a comparison of the computational cost with the spline-based method which is developed in \cite{deFrutos1}. Finally, Section \ref{Conclus} presents some concluding remarks.

All the algorithms presented in this work have been implemented in Matlab v2020b. All the numerical experiments have been performed in a personal computer with an Intel Core processor i7-8700K of 6 cores and 12 threads, with 3,70GHz(base)/4,70GHz(turbo) and 16Gb of RAM memory.

\section{A pollution differential game}\label{pollutiong}

The model is a J-player non cooperative differential game. Let $\Omega$ be a planar region with a given partition in $J$ subdomains such that
\begin{equation}
\overline{\Omega}=\bigcup^{J}_{j=1}\overline{\Omega}_j, \ \Omega_{i}\cap\Omega_{j}=\emptyset, \ i\neq j,
\end{equation}
where $\overline{\Omega}$ denotes the closure of $\Omega$.

Let $\partial_{ij}$ be the common boundary between $\Omega_i$ and $\Omega_j$, i.e.
\begin{equation}
\partial_{ij}=\partial\Omega_i\cap \partial\Omega_j=\overline{\Omega}_i\cap\overline{\Omega}_j, \ i\neq j.
\end{equation}

Each player $i$ controls just region $\Omega_i$ and he can choose the rate of pollutant emissions in that region. The objective of each player is to maximize his own payoff.

Let $u_i(\boldsymbol{x},t), \ i=1,...,J$ be the emission rate of subregion $i$, at time $t\geq0$ and at point $\boldsymbol{x}\in\Omega$. Function $P(\boldsymbol{x},t)$ denotes the stock of pollution defined $\forall \boldsymbol{x}\in \Omega$.

For scalar functions $f:\Omega\rightarrow \mathbb{R}$, symbol $\nabla f$ corresponds to the spatial gradient and, for vectorial functions $\boldsymbol{f}:\Omega\rightarrow \mathbb{R}^2$, symbol $\nabla\cdot f=\frac{\partial f_1}{\partial x}+\frac{\partial f_2}{\partial y}$ represents the divergence.

The main objective in \cite{deFrutos1} was to study the spatial relation between decision makers. We are going to stick at the simplest model (no wind pollution transport, no non-linear reaction terms, simplest discrete-space model version...). More complex models, which might require further numerical treatment, will be considered in future research (see Section \ref{Conclus}).

The following parabolic partial differential equation gives the spatio-temporal dynamics of the stock of pollution:
\begin{equation}\label{systemec}
\begin{aligned}
& \frac{\partial P}{\partial t}=\nabla \cdot (k\nabla P)-cP+F(\boldsymbol{u}), \quad \boldsymbol{x}\in\Omega, \\
& P(\boldsymbol{x},0)=P_0(\boldsymbol{x}), \quad \boldsymbol{x}\in\Omega,\\
&\alpha(\boldsymbol{x})P(\boldsymbol{x},t)+k(\boldsymbol{x})\nabla P^T(\boldsymbol{x},t)\boldsymbol{n}=\alpha(\boldsymbol{x})P_b(\boldsymbol{x},t), \quad \boldsymbol{x}\in \partial\Omega,
\end{aligned}
\end{equation}
where $\boldsymbol{u}=[u_1,...,u_J]^T$ is the vector of emission rates, $k=k(\boldsymbol{x})$ is a local diffusion coefficient, which is assumed to be a smooth function such that $k_m\leq k(\boldsymbol{x})\leq k_M, \ \forall \boldsymbol{x}\in\Omega$ and $0<k_m<k_M$ are given constants. This coefficient measures the velocity at which the stock of pollutant is diffused in a location $\boldsymbol{x}$. Term $cP=c(\boldsymbol{x},t)$ represents the natural decay of pollutant.

It is assumed that only agent $j$ emits in subregion $\Omega_j, \ j=1,...,J$ and that each $\boldsymbol{x}\in\Omega$ belongs to just one region. Therefore, the source term can be written as:
\begin{equation}
F(\boldsymbol{u}(\boldsymbol{x},t))=\sum^J_{j=1}F_j(u_j(\boldsymbol{x},t))\boldsymbol{1}_{\Omega_j}(\boldsymbol{x}),
\end{equation}
where $F_j, \ j=1,...,J$ are smooth functions and $\boldsymbol{1}_{\Omega_j}$ is the characteristic function of $\Omega_j$. By the hypothesis of the model, we have that $F(\boldsymbol{u}(\boldsymbol{x},t))=F_j(u_j(\boldsymbol{x},t))$ if $\boldsymbol{x}\in\Omega_j$.

Concerning boundary condition, $\alpha(\boldsymbol{x})$ is a non-negative smooth function that appears due to Newton's law of diffusion on the boundary of $\Omega$.

The objective of player $i, \ i=1,...,J$ is to maximize his payoff
\begin{equation}
J_i(u_1,...,u_J,P_0)=\int_0^{+\infty}\int_{\Omega_{i}}e^{-\rho t}G_i(u_1,...,u_J,P)d\boldsymbol{x}dt,
\end{equation}
subject to the dynamics given by (\ref{systemec}). Parameter $\rho>0$ is a given time-discount rate. The instant welfare $G$ is given by a benefit from consumption minus the damage caused by the stock of pollutants.

Each region $i$ produces one consumption good, where the amount of production is controlled by player $i$, and such production produces emissions (pollution). Therefore, we can represent
\begin{equation*}
G_i(u_1,...,u_J,P)=(B_i(u_i)-D_i(P))\boldsymbol{1}_{\Omega_i}
\end{equation*}
where $B_i(u_i)$ corresponds to the instantaneous benefits from production and $D_i(P)$ to the environmental damage caused by the accumulated stock of pollution. $B_i$ and $D_i$ are assumed to be smooth functions and respectively concave and convex in their arguments.

Now we proceed to describe the discrete-space version of the model. We only sketch the main ideas and we refer to Appendix B,\cite{deFrutos1} for the details.

Functions $u_i, P_i$ are considered densities of emissions and pollution stocks along region $\Omega_i$. We define
\begin{equation}
p_i(t)=\frac{1}{m_i}\int_{\Omega_i}P(\boldsymbol{x},t)d\boldsymbol{x}, \quad \quad v_i(t)=\frac{1}{m_i}\int_{\Omega_i}u_i(\boldsymbol{x},t)d\boldsymbol{x}, \ i=1,...,J
\end{equation}
where $m_i=\int_{\Omega_i}d\boldsymbol{x}$.

Under a linear-cuadratic specification and an infinite-time horizon
\begin{equation}
\begin{aligned}
& F_i(v_1,...,v_J):=\beta_iv_i, \quad \quad G_i(v_1,...,v_J,\boldsymbol{p}):=v_i\left(A_i-\frac{v_i}{2}\right)-\frac{\varphi_i}{2}p_i^2, \\
& \boldsymbol{p}=[p_1,...,p_J]^T, \quad v_i=v_i(\boldsymbol{p}), \quad m_i=m_j, \ \forall i,j=1,...,J,
\end{aligned}
\end{equation}
and some calculus, the objective of player $i$ is to maximize
\begin{equation}
J_i(v_1,...,v_J,\boldsymbol{p}_0)=\int_0^{+\infty}e^{-\rho t}\left(v_i\left(A_i-\frac{v_i}{2}\right)-\frac{\varphi_i}{2}p_i^2\right)dt,
\end{equation}
subject to the dynamics of the aggregated stock of pollution given by the set of ordinary differential equations
\begin{equation}
m_i\frac{dp_i}{dt}=\sum^{J}_{\substack{j=0 \\ j\neq i}}k_{ij}(p_i-p_j)-m_ic_ip_i+m_iF(v_i), \  i=1,...,J
\end{equation}
supplemented with a given initial state of pollution $\boldsymbol{p}^{0}=\left[p_1^0,...,p_j^{0}\right]^T$.

\subsection{A Chebyshev-based numerical method}\label{SecNummeth}

Let $h>0$ be a positive parameter, $t_n=nh$ the discrete times defined for all positive integers $n$ and $\delta_h=1-\rho h$ the discrete discount factor.

We denote by $\bar{u}_i, \ i=1,...,J$ a sequence of real numbers $\bar{u}_i=\{u_{i,n}\}_{n=0}^{\infty}$ and $\mathcal{U}$ denotes the set of real sequences $\bar{v}$ with $v_n\geq0, \ \forall n\in\mathbb{N}$.

For $\boldsymbol{p}=[p_1,...,p_J]^T\in\mathbb{R}^J$ and $\boldsymbol{u}=[u_1,...,u_J]^T\in\mathbb{R}^J, \ u_i\geq0, \ i=1,...,J$, we define
\begin{equation}
g_i(\boldsymbol{p},\boldsymbol{u})=\sum^{J}_{\substack{j=0\\ j\neq i}}\frac{k_{ij}}{m_i}(p_i-p_j)-c_ip_i+F(u_i), \  i=1,...,J
\end{equation}
and we denote $\boldsymbol{g}(\boldsymbol{p},\boldsymbol{u})=[g_1(\boldsymbol{p},\boldsymbol{u}),...,g_J(\boldsymbol{p},\boldsymbol{u})]^{T}$.

In the time-discrete infinite horizon game, each player $i=1,...,J$ wants to maximize
\begin{equation}
W_i(\bar{u}_i,\boldsymbol{p}_0)=h\sum^{\infty}_{n=1}\delta^n_{h}G_i(u_{i,n},p_{i,n}), \ \bar{u}_i\in\mathcal{U},
\end{equation}
subject to
\begin{equation}\label{ecgtimedisc}
\boldsymbol{p}_{n+1}=\boldsymbol{p}_{n}+h\boldsymbol{g}(\boldsymbol{p}_{n},\boldsymbol{u}_{n}), \ n\geq 0
\end{equation}
where $\boldsymbol{p}_{n}=[p_{1,n},...,p_{J,n}]^T, \ \boldsymbol{u}_{n}=[u_{1,n},...,u_{J,n}]^T$ and $\boldsymbol{p}_{0}$ is a given initial state.

The time-discrete value function $V_{h,i}(\boldsymbol{p}), \ i=1,...,J$ is obtained solving Bellman's equation
\begin{equation}
V_{h,i}(\boldsymbol{p})=\underset{u_i\geq0}{\max}\left\{hG_i(p_i,u_i)+\delta_hV_{h,i}\left(\boldsymbol{p}+h\boldsymbol{g}(\boldsymbol{p},[u_i,\boldsymbol{u}^{*}_{-i}])\right)\right\}
\end{equation}
where for $i=1,...,J$
\begin{equation}
\boldsymbol{u}^{*}_{i}=\underset{u_i\geq 0}{\text{argmax}}\left\{hG_i(p_i,u_i)+\delta_hV_{h,i}\left(\boldsymbol{p}+h\boldsymbol{g}(\boldsymbol{p},[u_i,\boldsymbol{u}^{*}_{-i}])\right)\right\}
\end{equation}
and where, from now on, we employ the notation
\begin{equation}
[u_i,\boldsymbol{v}_{-i}]=[v_1,...,v_{i-1},u_i,v_{i+1},...,v_{J}]^T, \ u_i\in \mathbb{R}, \ \boldsymbol{v}\in\mathbb{R}^J.
\end{equation}

We now present the main steps of a generalized collocation Chebyshev-based method. A review of Chebyshev interpolation and an implementation  is presented in Section 3.

\

\noindent\underline{Step 0:} Offline Computation

We define $N_p=(N^p_1,...,N^p_J)\in\mathbb{N}^J$ and $N_u=(N^u_1,...,N^u_J)\in\mathbb{N}^J$, two $J$-dimensional vectors such that $N^p_i,\ N^u_i>0, \ i=1,...,J$.

With these $J$-dimensional vectors, we build two adecuate sets of collocation points $\boldsymbol{P}\subset \mathbb{R}^{J}$, $\boldsymbol{U}\subset\mathbb{R}^{J}$ (detailed in Section  \ref{Jmulteval}).

Let $N_{\boldsymbol{P}}=\left|\boldsymbol{P}\right|$ and $\boldsymbol{P}=\left\{\bar{\boldsymbol{p}}_j\in\mathbb{R}^J, \ j=1,...,N_{\boldsymbol{P}}\right\}$.

For each player $i=1,...,J$, we compute a Chebyshev interpolation polynomial in the control variables for every collocation node in the state variables, i.e. we compute
\begin{equation*}
g^i_{\bar{\boldsymbol{p}}_j}({\boldsymbol{u}}), \quad j=1,...,N_{\boldsymbol{P}},
\end{equation*}
which are $N_{\boldsymbol{P}}$ different interpolation polynomials in $\boldsymbol{u}$, such that $\forall j=1,...,N_{\boldsymbol{P}}$ it holds
\begin{equation*}
g^{i}_{\bar{\boldsymbol{p}}_j} \left(\bar{\boldsymbol{u}}\right) =g_i\left(\bar{\boldsymbol{p}}_j,\left[\bar{u}_{i},\bar{\boldsymbol{u}}_{-i}\right]\right), \quad \forall  \bar{\boldsymbol{u}}\in\boldsymbol{U}
\end{equation*}

We denote $\boldsymbol{g}_{\bar{\boldsymbol{p}}_j}({\boldsymbol{u}})=\left[g^{1}_{\bar{\boldsymbol{p}}_j}\left({\boldsymbol{u}}\right),g^{2}_{\bar{\boldsymbol{p}}_j}\left({\boldsymbol{u}}\right),...,g^{J}_{\bar{\boldsymbol{p}}_j}\left({\boldsymbol{u}}\right)\right], \ j=1,2,...,N_{\boldsymbol{P}}$.

We compute some localization indexes (detailed in Section \ref{Jmulteval}).

We set $r=0$ and a small time step $h\in\mathbb{R}^{+}$.

For each player $i=1,...,J$, we initialize the iteration with some given $V^{N_p,[0]}_{h,i}\left(\bar{\boldsymbol{p}}_j\right)$ and $\boldsymbol{u}^{[0]}\left(\bar{\boldsymbol{p}}_j\right), \ j=1,2,...,N_{\boldsymbol{P}}$.

For each player $i=1,...,J$, we compute the Chebyshev interpolation polynomial $V^{N_p,[0]}_{h,i}(\boldsymbol{p})$ which interpolates $V^{N_p,[0]}_{h,i}\left(\bar{\boldsymbol{p}}_j\right), \ j=1,2,...,N_{\boldsymbol{P}}$.

\

\noindent\underline{Step 1:}

For each player $i=1,...,J$ and each $\bar{\boldsymbol{p}}_j, \ j=1,2,...,N_{\boldsymbol{P}}$, we compute the $J$-dimensional and one variable polynomial
\begin{equation*}
\mathcal{G}^i_{\bar{\boldsymbol{p}}_j}(u)=\left.\boldsymbol{g}_{\bar{\boldsymbol{p}}_j}(\boldsymbol{u})\right|_{\boldsymbol{u}^{[r]}_{-i}(\bar{\boldsymbol{p}}_j)}, \ j=1,2,...,N_{\boldsymbol{P}}
\end{equation*}

\noindent\underline{Step 2:}

For each player $i=1,...,J$ and each $\bar{\boldsymbol{p}}_j, \ j=1,2,...,N_{\boldsymbol{P}}$, we compute the one variable polynomial
\begin{equation*}
\mathcal{V}^{N_p,[r]}_{h,i,\bar{\boldsymbol{p}}_j}(u)=V^{N_p,[r]}_{h,i}\left({\bar{\boldsymbol{p}}_j+h\mathcal{G}^i_{\bar{\boldsymbol{p}}_j}(u)}\right).
\end{equation*}

\noindent\underline{Step 3:}

For each player $i=1,...,J$, we find the strategy at each state node $\bar{\boldsymbol{p}}_j, \ j=1,2,...,N_{\boldsymbol{P}}$ which maximizes the objective function, i.e.

\begin{equation*}
u^{[r+1]}_{i}\left(\bar{\boldsymbol{p}}_j\right)=\underset{u\geq 0}{\text{argmax}}\left\{\mathcal{V}^{N_p,[r]}_{h,i,\bar{\boldsymbol{p}}_j}(u)\right\}.
\end{equation*}

\noindent\underline{Step 4:}

For each player $i=1,...,J$, we define $V^{N_p,[r+1]}_{h,i}(\boldsymbol{p})$ as the Chebyshev interpolation polynomial which interpolates $\mathcal{V}^{N_p,[r]}_{h,i,\bar{\boldsymbol{p}}_j}\left(u^{[r+1]}_{i}\left(\bar{\boldsymbol{p}}_j\right)\right), \ j=1,2,...,N_{\boldsymbol{P}}$.

If we are not below the prescribed tolerance,
\begin{equation*}
\left|V^{N_p,[r+1]}_{h,i}(\boldsymbol{p})-V^{N_p,[r]}_{h,i}(\boldsymbol{p})\right|<\text{TOL}, \quad i=1,...,J
\end{equation*}
we set $r=r+1$ and return to Step 1. Otherwise, we stop.

\

We point out that, in the particular pollution problem we are dealing with, $g^{i}_{\bar{\boldsymbol{p}}_j}\left({\boldsymbol{u}}\right)=g^{i}_{\bar{\boldsymbol{p}}_j}\left(u_i\right), \ i=1,...,J, \ \forall \bar{\boldsymbol{p}}_j\in\boldsymbol{P}$ is one dimensional, but we prefer to present a generalized algorithm in the case it was not.

\section{The Chebyshev interpolation}\label{cbmethod}

We now make first a brief review of multidimensional Chebyshev interpolation and comment how the different calculus involved in the previous algorithm can be efficiently performed.

We are going to employ the work presented in Section 2,\cite{deFrutos2}, where it is described how multidimensional Chebyshev polynomials can be efficiently computed, storaged and evaluated for several values in all the dimensions simultaneously.

Here, we only include the main definitions in \cite{deFrutos2} and the modifications needed to adapt the algorithm to the problem described in Section \ref{pollutiong}.

\subsection{A review of multidimensional Chebyshev interpolation}\label{Reviewchpol}

The Chebyshev polynomial of degree $n$ (see \cite{Rivlin}) is given by
\begin{equation*}
T_n(x)=\cos\left(n \arccos (x) \right),
\end{equation*}
where $0\leq \arccos (x) \leq \pi$.

\

From now on, variable $x\in[-1,1]$ or $\textbf{x}=(x_1,...,x_n)\in[-1,1]^{n}$ for the $n$-dimensional case.

\

Let $N\in \mathbb{N}$. The $N+1$ Chebyshev nodes $\{\alpha^k\}_{k=0}^{N}$ in interval $[-1, \ 1]$ correspond to the extrema of $T_n(x)$ and they are given by:
\begin{equation*}
\alpha^k=\cos\left(\frac{\pi k}{N}\right), \quad k=0,1,...,N.
\end{equation*}

If the function $F(\tilde{x})$ that we want to interpolate is defined in interval $\tilde{x}\in[a,b]$, the Chebyshev nodes $\{\tilde{\alpha}^k\}_{k=0}^{N}$ in interval $[a,b]$ are computed with the $\{\alpha^k\}_{k=0}^{N}$ nodes in $[-1,1]$ and the change of variable given by formula
\begin{equation}
\tilde{x}=\frac{b-a}{2}x+\frac{b+a}{2}, \quad x\in[-1,1].
\end{equation}

\begin{definition}\label{Ch1defpolinter1variable}
Let $F(\tilde{x})$ be a continuous function defined in $\tilde{x} \in [a,b]$.

For $N\in\mathbb{N}$, let $I_N F(x)$ be the N degree interpolant of function  $F(\tilde{x})$ at the Chebyshev nodes, i.e. the polynomial which satisfies
\begin{equation*}
I_N F(\alpha^k)=F(\tilde{\alpha}^k), \quad k=0,1,...,N.
\end{equation*}

Polynomial $I_N F(x)$ can be expressed as
\begin{equation}
I_N F(x)=\sum^{N}_{l=0}\hat{p}_{l}T_l(x), \quad x\in[-1, 1],
\end{equation}
where coefficients $\hat{p}_{l}$ are given by
\begin{equation}\label{forcoef}
\begin{aligned}
\hat{p}_l&=\frac{1}{N}\sum^{N}_{k=0}{}^{''}F(\tilde{\alpha}^k)T_l(\alpha^k), \quad \text{if } l\in\{0,N\}, \\
\hat{p}_l&=\frac{2}{N}\sum^{N}_{k=0}{}^{''}F(\tilde{\alpha}^k)T_l(\alpha^k), \quad \text{if } l\in\{1,2,...,N-1\},\\
\end{aligned}
\end{equation}
and the double prime indicates that we halve the first and last elements.
\end{definition}

Instead of using formula (\ref{forcoef}), we will employ an efficient FFT based algorithm which is presented in \cite{Canuto} or \cite{deFrutos2}. For the univariate case

\noindent{\textit{Algorithm C1v:}}

1. Define
\begin{equation*}
z=\left[F\left(\tilde{\alpha}^0\right),F\left(\tilde{\alpha}^1\right),...,F\left(\tilde{\alpha}^{N-1}\right),F\left(\tilde{\alpha}^N\right),F\left(\tilde{\alpha}^{N-1}\right),...,F\left(\tilde{\alpha}^1\right)\right]^T
\end{equation*}

2. Compute
\begin{equation*}
y=\frac{real(FFT(z))}{2N}
\end{equation*}

3. It holds that
\begin{equation*}
\left\{\begin{aligned}
\hat{p}_0 &=y(1), \\
\hat{p}_l &=y(l+1)+y(2N-(l-1)), \ \text{if } 0<l<N, \\
\hat{p}_N &=y(N)
\end{aligned}
\right.
\end{equation*}

We also mention the algorithm presented in \cite{Canuto} which allows to compute efficiently the derivative of a Chebyshev interpolation polynomial.

\begin{proposition}
If $F(\tilde{x})$ is a continuous function defined in $\tilde{x} \in [a,b]$ and
\begin{equation*}
I_{N}F(x)=\sum_{l=0}^{N}\hat{p}_lT_l(x), \ x\in[-1,1]
\end{equation*}
is its Chebyshev interpolation polynomial, it holds that
\begin{equation*}
\left(I_{N}F(x)\right)'=\frac{2}{b-a}\sum_{l=0}^{N-1}\hat{q}_lT_l(x)
\end{equation*}
where for $l=0,1,...,N-1$:
\begin{equation*}
\hat{q}_l=\frac{2}{c_l}+\sum_{\substack{\tiny j=l+1\\j+l \ odd \normalsize}}^Nj\hat{p}_j, \ \text{where} \ c_l=\left\{\begin{aligned}2, \ l=0,\\ 1, \ l\geq1.\end{aligned}\right.
\end{equation*}
\end{proposition}

Now we proceed to multidimensional interpolation.

\begin{definition}\label{Ch1defpolinternvariable}
Let $\tilde{\textbf{x}}=\left(\tilde{x}_1,\tilde{x}_2,...,\tilde{x}_n \right)$ and $\tilde{F}(\tilde{\textbf{x}})$ be a continuous function defined in $\tilde{x}_j\in[a_j, \ b_j], \quad j=1,2,...,n$.

For $\boldsymbol{N}=\{N_1,N_2,...,N_n\}\in\mathbb{N}^n$, we define
\begin{equation}
L^{\boldsymbol{N}}=\left\{\boldsymbol{l}=(l_1,l_2,...,l_n) \ / \ 0\leq l_j\leq N_j, \ l_j\in\mathbb{N}, \ j=1,2,...,n\right\}.
\end{equation}

For $j=1,2,...,n$, let $\left\{\alpha^k_j\right\}_{k=0}^{N_j}$ be the $N_j+1$ Chebyshev nodes in $[-1, 1]$ and $\left\{\tilde{\alpha}^k_j\right\}_{k=0}^{N_j}$ the corresponding $N_j+1$ Chebyshev nodes in $[a_j, b_j]$.

\

We use the notation $\tilde{\boldsymbol{\alpha}}^{\boldsymbol{l}}=\left(\tilde{\alpha}^{l_1}_1,\tilde{\alpha}^{l_2}_2,...,\tilde{\alpha}^{l_n}_n \right)$ and $\boldsymbol{\alpha}^{\boldsymbol{l}}= \left(\alpha^{l_1}_1, \ \alpha^{l_2}_2,...,\alpha^{l_n}_n \right)$.

\

Let $I_{\boldsymbol{N}} F(\textbf{x})$ be the n-dimensional interpolant of function $F(\tilde{\textbf{x}})$ at the Chebyshev nodes, i.e. the polynomial which satisfies
\begin{equation*}
I_{\boldsymbol{N}} F(\boldsymbol{\alpha}^{\boldsymbol{l}})=F(\tilde{\boldsymbol{\alpha}}^{\boldsymbol{l}}), \quad \boldsymbol{l}\in L^{\boldsymbol{N}}.
\end{equation*}

Polynomial $I_{\boldsymbol{N}} F(\textbf{x})$ can be expressed as
\begin{equation}
I_{\boldsymbol{N}} F(\textbf{x})=\sum_{\boldsymbol{l}\in L^{\boldsymbol{N}}}\hat{p}_{\boldsymbol{l}}T^{\boldsymbol{l}}(\textbf{x}), \quad \textbf{x}\in[-1, 1]^{n},
\end{equation}
where
\begin{equation*}
T^{\boldsymbol{l}}(\textbf{x})=T_{l_1}(x_1)T_{l_2}(x_2) ... T_{l_n}(x_n).
\end{equation*}
and the coefficients $\hat{p}_{\boldsymbol{l}}=\hat{p}_{(l_1,l_2,...,l_n)}\in \mathbb{R}$ can be computed with the $n$-dimensional version of the Algorithm C1v presented before.
\end{definition}

\noindent{\textit{Algorithm Cnv:}}

Let $\Gamma_{(N_1+1)\times...\times(N_n+1)}$ be a n-dimensional array such that
\begin{equation*}
\Gamma(l_1+1,l_2+1,...,l_n+1)=F(\tilde{\alpha}^{l_1}_1,\tilde{\alpha}^{l_2}_2,...,\tilde{\alpha}^{l_n}_n)
\end{equation*}

\noindent 1. $A_1=\Gamma$.
\noindent 2. For $i=1$ to $n$

\hspace{0.5 cm} 2.1. $\{m_1,m_2, ...,m_n\}=\dim(B_i)$.

\hspace{0.5 cm} 2.2. For $j_2=1$ to $m_2$, for $j_3=1$ to $m_3$, ...,  for $j_n=1$ to $m_n$
\begin{equation*}
B_i(:,j_2,j_3,...,j_n)=\text{Algorithm C1v}\left(A_i(:,j_2,j_3,...,j_n)\right).
\end{equation*}

\hspace{0.5 cm} 2.3. $A_{i+1}=\text{permute}(B_i),[2:n \ 1])$.

\

\noindent 3. $\hat{p}_{\boldsymbol{l}}=A_{n+1}(l_1+1,l_2+1,...,l_n+1)$.

\

We remark that the FFT routine in Matlab admits multidimensional evaluation, so step 2.2 can be efficiently computed without loops.

Therefore, the polynomial coefficients are stored in a $(N_1+1)\times...\times(N_n+1)$-dimensional array $A$, where
\begin{equation*}
A(l_1+1,l_2+1,...,l_n+1)=\hat{p}_{(l_1,l_2,...,l_n)}
\end{equation*}

\subsection{Evaluation of one $N_u$-dimensional polynomial} \label{ndimeval}

Suppose now that we have a Chebyshev interpolation polynomial $I_{N_u} g(\textbf{u})$, given by a $(N^u_1+1)\times...\times(N^u_n+1)$-dimensional array $A$ and we want to evaluate it in a set of points $\{b^1_j\}_{j=1}^{k_1}$ just in the first variable, i.e. we want to compute
\begin{equation*}
\left\{I_{\boldsymbol{N}} F\left(b^1_j,u_2,u_3,...,u_n\right)\right\}_{j=1}^{k_1}.
\end{equation*}

In (Section 2,\cite{deFrutos2}) it is described how $\{(T_{l_1}(b^1_1),...,T_{l_1}(b^1_{k_1}))\}_{l_1=0}^{N_1}$ can be efficiently evaluated and stored in a ($k_1,N_1+1$)-dimensional array B such that
\begin{equation*}
B(j,l)=T_{l}(b^1_j)
\end{equation*}

Afterwards, a standard matrix product has to be performed over all the other dimensions. We need to compute
\begin{equation*}
B\cdot A(:,i_2,...,i_n), \ i_s=1,...,N_s+1, \ s=2,...,N.
\end{equation*}

In the last version of Matlab, this can be efficiently performed with ``\textit{pagemtimes}'' function. We can define
\begin{equation*}
C=permute(pagemtimes(B,A),[2:N_n \ 1])
\end{equation*}
where the result is a $(N_2+1)\times...\times(N_n+1)\times k_1$ dimensional array. The permutation is needed in order to evaluate further dimensions.

Array $C$ corresponds to the coefficients of the interpolation polynomial $I_{\boldsymbol{N}} F(\textbf{x})$ evaluated in the points $\{b^1_j\}_{j=1}^{k_1}$, i.e.
\begin{equation*}
C(:,...,:,j)\sim I_{\boldsymbol{N}} F(b^1_j,u_2,...,u_n), \ j=1,...,k_1.
\end{equation*}

If we want now to evaluate the polynomial in a set of points $\{b^2_j\}_{j=1}^{k_2}$ in the second variable, another set of points in the third variable..., we would proceed iteratively obtaining, at the end, a $(k_1,...,k_n)$-dimensional array $D$ which contains the evaluation of the polynomial in every possible combination of the points of each variable, i.e.

\begin{equation*}
D(j_1,j_2,...,j_n)=I_{\boldsymbol{N}} F(b^1_{j_1},b^2_{j_2},...,b^n_{j_n}), \ j_s=1,...,k_s, \ s=1,...,n
\end{equation*}

\subsection{Evaluation of $N_{\boldsymbol{P}}$ different $N_u$-dimensional polynomials in different points} \label{multevalmulpol}

Suppose that we have $N_{\boldsymbol{P}}$ different multidimensional Chebyshev interpolation polynomials, where each one is given with a $N_u=(N^u_1+1,...,N^u_J+1)$-dimensional array $A_j, \ j=1,...,N_{\boldsymbol{P}}$ as shown in Subsection \ref{Reviewchpol}.

They can all be stored in a $(N^u_1+1,...,N^u_J+1,N_{\boldsymbol{P}})$-dimensional array $A_j$ where
\begin{equation*}
A(:,...,:,j)=A_j\sim I_{N_u} g_j(\boldsymbol{u}), \ j=1,...,N_{\boldsymbol{P}}
\end{equation*}
and $g_j(\boldsymbol{u}), \ j=1,...,N_{\boldsymbol{P}}$ is each of the functions that has been interpolated.

\

In order the employ the algorithm of Subsection \ref{ndimeval} efficiently in our pollution problem, a small modification has to be done.

\

Suppose that we want to evaluate each polynomial in a different point in the first variable, i.e., given $\{b^1_{j}\}_{j=1}^{N_{\boldsymbol{P}}}$ we have to compute
\begin{equation*}
\left\{I_{N_u} g_{j}(b^1_j,u_2,u_3,...,u_n)\right\}_{j=1}^{N_{\boldsymbol{P}}},
\end{equation*}

We remark that in Subsection \ref{ndimeval} we wanted to evaluate (in the first variable) one polynomial in a set of $k_1$ different points. Here we want to evaluate each polynomial $g_j(b^1_j,u_2,u_3,...,u_n)$ in a specific point $b_j, \ j=1,...,N_{\boldsymbol{P}}$.

\

We build a 2-dimensional array $B$ as defined in Subsection \ref{Reviewchpol} such that $B(j,l)=T_{l}(b^1_j)$, and we define the following location index
\begin{equation*}
\begin{aligned}
&aux1=[1:N_{\boldsymbol{P}}:N_{\boldsymbol{P}}(N^u_2+1)...(N^u_J+1)]\\
&locind_{1}=aux1\\
&for \ l=2:N_{\boldsymbol{P}} \\
& \hspace{1 cm} aux2=N_{\boldsymbol{P}}(N^u_2+1)...(N^u_J+1)(l-1)+(l-1) \\
& \hspace{1 cm} locind_{1}=[locind_{1} \ \ (aux1+aux2)]\\
&end
\end{aligned}
\end{equation*}

The evaluation
\begin{equation*}
\begin{aligned}
& C=permute\left(pagemtimes\left(B,A\right),[2:J \ 1]\right) \\
& D=reshape(C(locind_{1}),[N^{u}_2 \ N^u_3 \ ... \ N^u_J \ N_{\boldsymbol{P}}])
\end{aligned}
\end{equation*}
gives a $(N^u_2+1,N^u_3+1,...,N^u_n+1,N_{\boldsymbol{P}})$-dimensional array $D$ where
\begin{equation*}
D(:,...,:,j)\sim I_{\boldsymbol{N}} g_j(b^1_j,u_2,...,u_J), \ i=1,...,J
\end{equation*}

In a similar way, a location index $locind_{2}$ can be computed to compute $\left\{I_{N_u} g_j(b^1_j,b^2_j,u_3...,u_J)\right\}_{j=1}^{N_{\boldsymbol{P}}}, \ i=1,...,J$ for a second set of points $\{b^2_{j}\}_{j=1}^{N_{\boldsymbol{P}}}$. And so on for evaluating the rest of the dimensions.

\subsection{Implementation of the algorithm} \label{Jmulteval}

\noindent \underline{\textit{Step 0:}} Offline computations

\

Suppose that the $J$ players are indexed by $i=1,...,J$.

Let $N_p=(N^p_1,...,N^p_J)\in\mathbb{N}^J$ and $N_u=(N^u_1,...,N^u_J)\in\mathbb{N}^J$ be two $J$-dimensional vectors such that $N^p_i,\ N^u_i>0, \ i=1,...,J$.

Vectors $N_p$ and $N_u$ will be respectively employed to define the discretization in the state space and in the control space.

\

Let us introduce two positive parameter $P_M,U_M>0$ big enough and consider intervals $[0,P_M]$ and $[0,U_M]$. For each player $i$, the Chebyshev nodes $\left\{\tilde{p}^i_j\right\}_{j=0}^{N^p_i}$ and $\left\{\tilde{u}^i_j\right\}_{j=0}^{N^u_i}$ are given by
\begin{equation*}
\begin{aligned}
\tilde{p}^i_j&=\frac{1}{2}\left[\cos\left(\frac{\pi j}{N^p_i}\right)(P_M-0)+(P_M+0)\right], \quad j=0,1,...,N^p_i, \\
\tilde{u}^i_j&=\frac{1}{2}\left[\cos\left(\frac{\pi j}{N^u_i}\right)(U_M-0)+(U_M+0)\right], \quad j=0,1,...,N^u_i,
\end{aligned}
\end{equation*}

We consider the J-intervals
\begin{equation*}
\begin{aligned}
\tilde{I}_p&=[0,P_M]\times...\times[0,P_M]\subset \mathbb{R}^{J} \\
\tilde{I}_u&=[0,U_M]\times...\times[0,U_M]\subset \mathbb{R}^{J}
\end{aligned}
\end{equation*}
where we will numerically solve the pollution game. We define the sets of collocation points
\begin{equation*}
\begin{aligned}
\tilde{\boldsymbol{P}}=\left\{\left(\tilde{p}^1_{j_1},\tilde{p}^2_{j_2},...,\tilde{p}^J_{j_n}\right), \ j_i=0,1,...,N^p_i, \ i=1,...,J\right\} \\
\tilde{\boldsymbol{U}}=\left\{\left(\tilde{u}^1_{j_1},\tilde{u}^2_{j_2},...,\tilde{u}^J_{j_n}\right), \ j_i=0,1,...,N^u_i, \ i=1,...,J\right\}
\end{aligned}
\end{equation*}

For simplicity in the notation we believe that, prior to initialize the algorithm, it is better to perform the corresponding changes of variables to $[-1,1]$ (as seen in Subsection \ref{Reviewchpol}).

\

Therefore, we will work directly with the J-intervals $I_p=I_u=[-1,1]^{J}$ and the corresponding sets of chebyshev collocation points
\begin{equation*}
\begin{aligned}
\boldsymbol{P}=\left\{\left(p^1_{j_1},p^2_{j_2},...,p^J_{j_n}\right), \ j_i=0,1,...,N^p_i, \ i=1,...,J\right\} \\
\boldsymbol{U}=\left\{\left(u^1_{j_1},u^2_{j_2},...,u^J_{j_n}\right), \ j_i=0,1,...,N^u_i, \ i=1,...,J\right\}
\end{aligned}
\end{equation*}
defined in $[-1,1]^J$. Once the algorithm is finished, we move back to the original intervals $\tilde{I}_p$ and $\tilde{I}_u$.

\

Therefore $N_{\boldsymbol{P}}=\left|\boldsymbol{P}\right|=\prod^J_{i=1}(N^p_i+1)$ and  $\boldsymbol{P}=\left\{\bar{\boldsymbol{p}}_j, \ j=1,...,N_{\boldsymbol{P}}\right\}$.

\

For any player $i\in\{1,...,J\}$, we need to compute $N_{\boldsymbol{P}}$ different interpolation polynomials of $\left\{g^{i}_{\bar{\boldsymbol{p}}}(\boldsymbol{u}), \ \bar{\boldsymbol{p}}\in\boldsymbol{P} \right\}$ such that $\forall \bar{\boldsymbol{p}}\in\boldsymbol{P}$, it holds
\begin{equation*}
g^{i}_{\bar{\boldsymbol{p}}} \left(\bar{\boldsymbol{u}}\right) =g_i\left(\bar{\boldsymbol{p}},\left[\bar{u}_{i},\bar{\boldsymbol{u}}_{-i}\right]\right), \quad \forall  \bar{\boldsymbol{u}}\in\boldsymbol{U}
\end{equation*}

We remark that these polynomials have to be computed just once and this can be efficiently done with \textit{Algorithm Cnv} as seen in Subsection \ref{Reviewchpol}. The polynomials will be $(N^u_1+1,...,N^u_J+1)$-dimensional and, for the rest of the algorithm, we identify for any player $i\in\{1,...,J\}$
\begin{equation*}
g^{i}_{\bar{\boldsymbol{p}}_j}(\boldsymbol{u}) \sim \left\{I_{N_u} g^{i}_j(\textbf{u})\right\}_{j=1}^{N_{\boldsymbol{P}}}, \ j=1,...,N_{\boldsymbol{P}}.
\end{equation*}

In the iterative algorithm, at any iteration $r$ and for any player $i\in\{1,...,J\}$, we will need to evaluate these polynomials in\small
\begin{equation*}
\left\{I_{N_u} g^{i}_{j}(u^{[r]}_1(\bar{\boldsymbol{p}}_j),u^{[r]}_2(\bar{\boldsymbol{p}}_j),...,u^{[r]}_{i-1}(\bar{\boldsymbol{p}}_j),u^{i}_{k},u^{[r]}_{i+1}(\bar{\boldsymbol{p}}_j),...,u^{[r]}_{J}(\bar{\boldsymbol{p}}_j))\right\}_{k=0}^{N^u_i}, \ j=1,...,N_{\boldsymbol{P}}
\end{equation*}\normalsize
where we recall that $\{u^{i}_{k}, \ k=0,...,N^u_{i}\}$ are the control Chebyshev nodes of player $i$.

\

Therefore, we can build a set of location indexes $locind_j, \ j=1,...,J$ which allow to perform such computation efficiently as shown in Subsection \ref{multevalmulpol}.

We remark that this location indexes have to be computed just once and can be employed in any iteration $[r]$ of the algorithm.

\

We initialize with some given $V^{N_p,[0]}_{h,i}(\bar{\boldsymbol{p}})$ and $\boldsymbol{u}^{[0]}(\bar{\boldsymbol{p}}_j), \ \bar{\boldsymbol{p}}\in\boldsymbol{P}$.

For each player $i=1,...,J$, we compute the Chebyshev interpolation polynomial $V^{N_p,[0]}_{h,i}(\boldsymbol{p})$, which interpolates $V^{N_p,[0]}_{h,i}(\bar{\boldsymbol{p}}), \ {\bar{\boldsymbol{p}}}\in\bar{\boldsymbol{P}}$ with \textit{Algorithm CnV}.

\

\noindent \underline{\textit{Step 1} and \textit{Step 2}:}

For every player $i\in\{1,..,J\}$ we compute $\left\{g^{i}_{\bar{\boldsymbol{p}}_j}\left(u^{i}_{k},\boldsymbol{u}^{[r]}_{-i}(\bar{\boldsymbol{p}}_j)\right)\right\}_{k=0}^{N^u_i}$, i.e.
\begin{equation}\small
\left\{I_{N_u} g^{i}_{j}(u^{[r]}_1(\bar{\boldsymbol{p}}_j),u^{[r]}_2(\bar{\boldsymbol{p}}_j),...,u^{[r]}_{i-1}(\bar{\boldsymbol{p}}_j),u^{i}_{k},u^{[r]}_{i+1}(\bar{\boldsymbol{p}}_j),...,u^{[r]}_{J}(\bar{\boldsymbol{p}}_j))\right\}_{k=0}^{N^u_i}, \ j=1,...,N_{\boldsymbol{P}}
\end{equation}\normalsize
with the technique described in Subsection \ref{multevalmulpol} and the location indexes precomputed in \textit{Step 0}.

\

We define
\begin{equation*}
\left\{\mathcal{G}^i_{\bar{\boldsymbol{p}}_j}\left(u^{i}_{k}\right)\right\}_{k=0}^{N^u_i}=\left\{\boldsymbol{g}_{\bar{\boldsymbol{p}}_j}\left(u^{i}_{k},\boldsymbol{u}^{[r]}_{-i}(\bar{\boldsymbol{p}}_j)\right)\right\}_{k=0}^{N^u_i}, \ j=1,2,...,N_{\boldsymbol{P}}
\end{equation*}
where we recall $\boldsymbol{g}_{\bar{\boldsymbol{p}}_j}({\boldsymbol{u}})=\left[g^{1}_{\bar{\boldsymbol{p}}_j}\left({\boldsymbol{u}}\right),g^{2}_{\bar{\boldsymbol{p}}_j}\left({\boldsymbol{u}}\right),...,g^{J}_{\bar{\boldsymbol{p}}_j}\left({\boldsymbol{u}}\right)\right], \ j=1,2,...,N_{\boldsymbol{P}}$.

\

We point out that, in practice, it is not necessary to build the interpolation polynomial of  $\left\{\mathcal{G}^i_{\bar{\boldsymbol{p}}_j}\left(u^{i}_{k}\right)\right\}_{k=0}^{N^u_i}$. For every $\bar{\boldsymbol{p}}\in\boldsymbol{P}$, in order to build $\mathcal{V}^{N_p,[r]}_{h,i_0,\bar{\boldsymbol{p}}}(u)$ we just compute
\begin{equation*}
V^{N_p,[r]}_{h,i_0}\left({\bar{\boldsymbol{p}}+h\mathcal{G}^{i}_{\bar{\boldsymbol{p}}}\left(u^{i_0}_{k}\right)}\right), \ k=0,1,...,N^u_{i_0}
\end{equation*}
and then apply \textit{Algorithm C1v} to the results obtained.

\

We want to remark that, working with arrays, all the operations can be implemented simultaneously for every $\bar{\boldsymbol{p}}\in\boldsymbol{P}$.

\

\noindent \underline{\textit{Step 4:}}

For any player $i_0\in\{1,..,J\}$, in order to compute
\begin{equation*}
u^{[r+1]}_{i}\left(\bar{\boldsymbol{p}}\right)=\underset{u\geq 0}{\text{argmax}}\left\{\mathcal{V}^{N_p,[r]}_{h,i,\bar{\boldsymbol{p}}}(u)\right\}, \ \bar{\boldsymbol{p}}\in\boldsymbol{P}.
\end{equation*}
we recommend to employ Newton algorithm for two reasons.

It is straightforward to implement Newton algorithm for all $\bar{\boldsymbol{p}}\in\boldsymbol{P}$ at the same time and the derivative of a Chebyshev interpolation polynomial can be efficiently obtained employing the algorithm presented in Subsection \ref{Reviewchpol}.

\subsection{Parallelization}\label{paralproc}

Since the evaluation over the $N_{\boldsymbol{P}}$ different state nodes is independent, the multidimensional arrays involved in the numerical algorithm described in Subsection \ref{multevalmulpol} can be split in smaller packages to different cores (computer processing units).

In our case, let $N_b$ and $N_f$ be two natural numbers such that $N_fN_b=N_{\boldsymbol{P}}$. For any array $A(:,...,:,1...N_{\boldsymbol{P}})$, employing \textit{reshape} function, we can redefine the array
\begin{equation*}
A=reshape(A,[N_1,...,N_J,N_f,N_b])
\end{equation*}

For $k=1,...,N_b$, we define $A'_k(:,...,:,1...N_f):=A(:,...,:,1...N_f,k)$.

\

The calculus involved in the numerical algorithm, for example the computation of $\left\{g^{i}_{\bar{\boldsymbol{p}}_j}\left(u^{i}_{k},\boldsymbol{u}^{[r]}_{-i}(\bar{\boldsymbol{p}}_j)\right)\right\}_{k=0}^{N^u_i}$ in Step 1, can be done independently in different cores employing Matlab \textit{parfor} and arrays $A'_k(:,...,:,1...N_f), \ k=1,...,N_b$. The information can be reassembled when needed.

The precomputation of localization indexes has also to be adapted to the smaller arrays that we have just defined, but this is something straightforward to do.

\

This parallelization procedure can also be applied working with just one core. If array $A(:,...,:,1...N_{\boldsymbol{P}})$ is very big, it can be splitted in smaller arrays as we have just described and solved with a standard \textit{for} loop.

\

The optimal (computing time) values for $N_b$ and $N_f$ depend on the values of $N_p$ and $N_u$, but probably they also depend on the number of cores and the kind of processors of the computer employed.

\begin{figure}[h]
\centering
\includegraphics[width=7cm,height=5 cm]{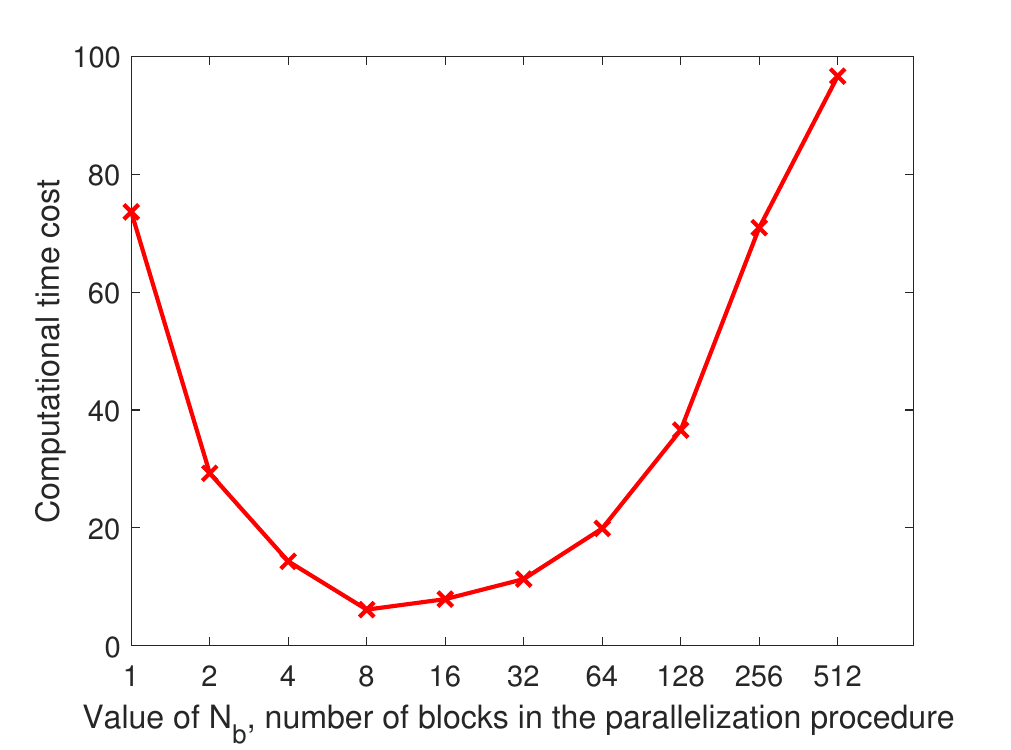}
\caption[Computational time cost for different size blocks in the parallelization.]{\label{exper2arallel} Computational time cost for different size blocks in the parallelization.}
\end{figure}

For example, with the computer that we employed in our experiments, we run a 3 players game with $N^p_i=7$ ($N_{\boldsymbol{P}}=512$). We computed the computational time cost of the numerical solution for smaller arrays given by $N_b=1,2,...,2^{9}$. The results are represented in Figure \ref{exper2arallel}.

This experiment shows that it was neither optimal to compute each state node in a different core (fully parallelization) nor to compute all the nodes at the same time in just one core (without parallelization and fully tensorized). The optimal computational time cost was ``half way'' between the size of the arrays involved and the number of blocks (which depends on the size of the arrays). Similar results were obtained when the game was played with different amounts of players.

\section{Numerical Results}\label{numres}

We now repeat some of the numerical experiments performed in \cite{deFrutos1}. We compare the spline method employed in that paper with the Chebyshev method that we have described.

When the pollution game is played by 2 players we have explicit solutions, so an error vs computational time cost analysis can be performed. For the case of 3 or more players, we lack of an explicit solution. We have obtained the same qualitative solutions as in \cite{deFrutos1}, but just a comparison of the computational time cost has been done.

\

Concerning the parallelization procedure, once we have the number of state nodes $N_{\boldsymbol{P}}$, let $\{M_1,...,M_{\sigma_0(N_{\boldsymbol{P}})}\}$ be all the natural dividers of $N_{\boldsymbol{P}}$.

For each numerical experiment, all the possible combinations for $N_f=M_i$ and $N_b=M_j$ such that $N_fN_b=N_{\Psi}$ have been tested. We point out  that for all the experiments,
\begin{itemize}
\item Case $N_f=N_{\boldsymbol{P}}$, $N_{b}=1$ (without parallelization and fully tensorized) is suboptimal.
\item Case $N_f=1$, $N_{b}=N_{\boldsymbol{P}}$  (fully parallel) is suboptimal.
\end{itemize}

The optimal computational time cost is always attained at some value $N_f=M_i$, $M_i\neq\{1,N_{\boldsymbol{P}}\}$.

\subsection{2 players}

We repeat Example 1 in \cite{deFrutos1}. Let
\begin{equation*}
\beta_i=1, \quad \varphi_i=1, \quad A_i=0.5, \quad c_i=0.5, \quad i=1,2, \quad K=[k_{ij}]=\left[\begin{matrix} -1 & 1 \\ 1 & -1\end{matrix}\right]
\end{equation*}

\

The spatial configuration described by $K$ means that players 1 and 2 share a common boundary and are isolated from outside.

We have computed the numerical solution for
\begin{itemize}
\item $h\in\{10^{-2},10^{-3},10^{-4},10^{-5}\}$,
\item $\text{TOL}\in\{10^{-2},10^{-3},10^{-4},10^{-5},10^{-6}\}$,
\item $N_i^p\in \{2,4,8\}, i=1,2.$
\end{itemize}

\begin{figure}[h]
\centering
\includegraphics[width=13cm,height=5 cm]{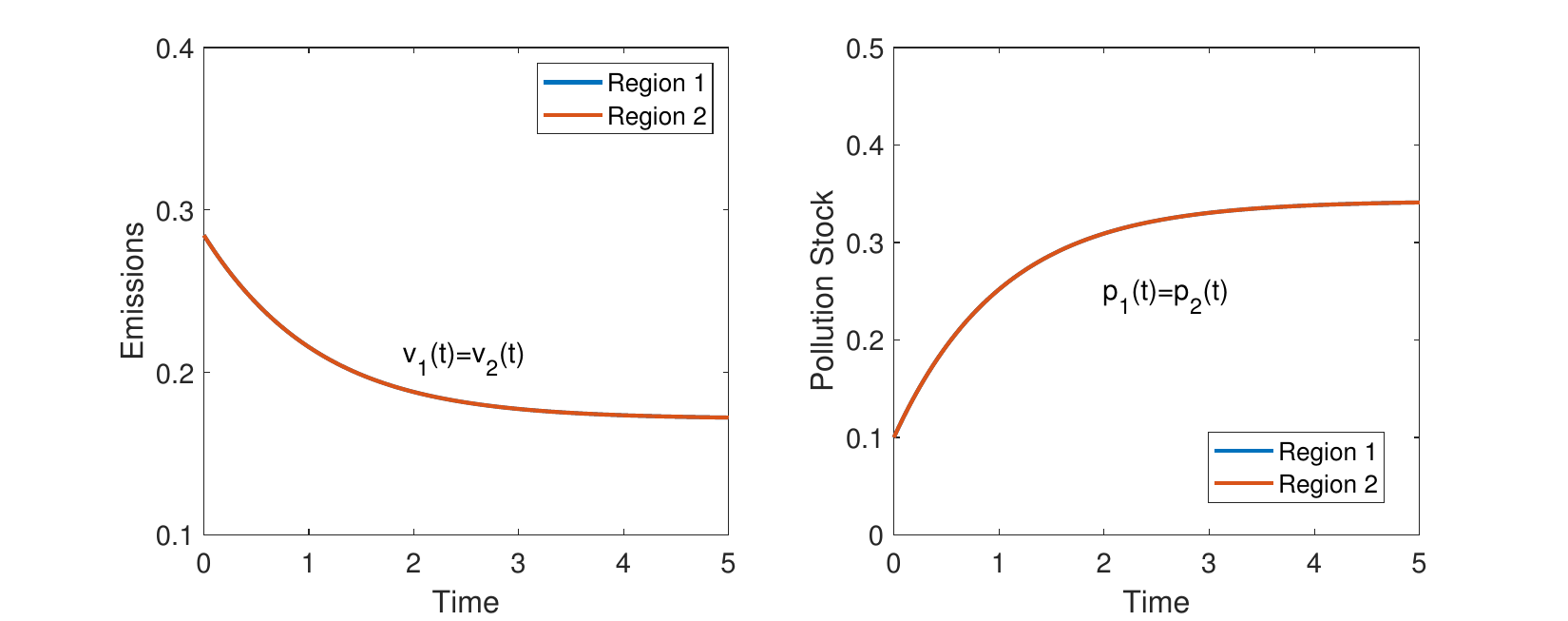}
\caption[Numerical emission (left) and pollution stock (right) time-paths along the equilibrium strategy obtained with the Chebyshev method.]{\label{exper1numsol} Numerical emission (left) and pollution stock (right) time-paths along the equilibrium strategy obtained with the Chebyshev method.}
\end{figure}

Under the spatial configuration defined, both players are symmetric, therefore the solutions of both players must coincide. In Figure \ref{exper1numsol} we represent the emission (left) and pollution (right) time paths obtained with the Chebyshev numerical method.

In order to analyse the performance, we study the numerical solution for the different values of $N^p_i, \ TOL$ and $h$.

For the 2 players case, we have explicit solutions (see \cite{deFrutos1}), so we can compute the exact optimal policy $u(x)$. For each experiment, we define the mean square error of the numerical solution by
\begin{equation*}
\text{error}=\frac{1}{N_{\Psi}}\sqrt{\sum_{x\in\Psi}\left(u^{*}(x)-u(x)\right)^{2}}
\end{equation*}
where $u^{*}$ is the numerical optimal policy obtained at the last iteration of the method in each experiment.

With the errors computed for all the experiments, we can plot the numerical error vs the computational time cost of each experiment and then retain the lower convex envolvent of the resulting cloud of points.

\begin{figure}[h]
\centering
\includegraphics[width=7cm,height=5 cm]{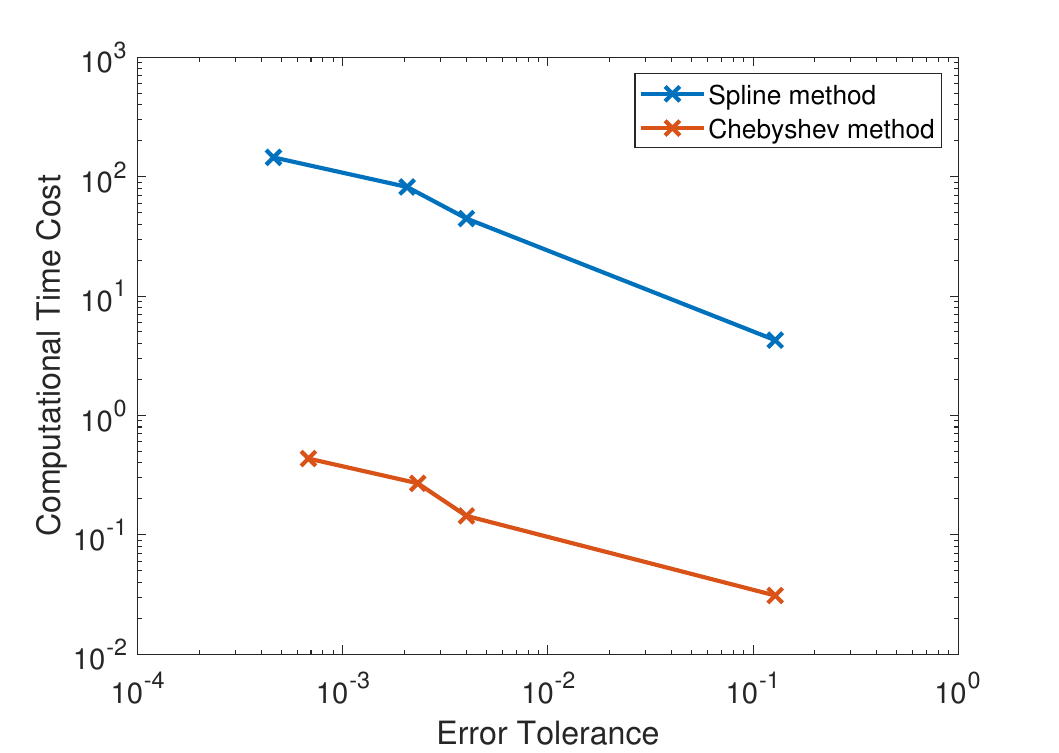}
\caption[Lower convex envolvent of the numerical error vs the computational time cost for the spline(blue) and Chebyshev(red) methods.]{\label{exper1perfor} Lower convex envolvent of the numerical error vs the computational time cost for the spline(blue) and Chebyshev(red) methods.}
\end{figure}

The lower convex envolvent informs, for a desired error tolerance, the minimum time required to attain that error. The analysis is represented in Figure \ref{exper1perfor}, for the spline(blue) and Chebyshev(red) methods.

The results in Figure \ref{exper1perfor} show that the Chebyshev method is much more efficient that the spline method. In average, for 2 players and  a similar prescribed error tolerance, the Chebyshev method requires 1/271 of the time of the spline method. The nodes of the lower convex with the biggest errors (the two situated at the right side) correspond to $N^p_i=2$, the next node to $N^p_i=4$ and the node with the smallest error (left side) corresponds to $N^p_i=8$.

It is interesting that both methods present the same error behaviour (the slopes of the lower convex envolvents are similar), since Chebyshev interpolation usually has a better error convergence than spline interpolation. This is probably due to the fact that the objective function has a linear-cuadratic specification and, therefore, both methods have similar error behaviour. It is possible that with non-polynomial objective specifications Chebyshev method could also present a better behaviour.

\subsection{3 players}

We now repeat Example 3 in \cite{deFrutos1}. The parameter values remain the same as in the previous experiment and the spatial configuration is given by

\begin{equation*}
K=[k_ij]=\left[\begin{matrix} -1 & 1 & 0 \\ 1 & -2 & 0 \\ 0 & 1 & -1\end{matrix}\right]
\end{equation*}

This configuration means that Player 2 shares a boundary with both Players 1 and 3, Players 1 and 3 have no common boundary and all the countries are isolated from outside. Under this configuration, Players 1 and 3 are symmetric, so their strategies should coincide.

\begin{figure}[h]
\centering
\includegraphics[width=13cm,height=5 cm]{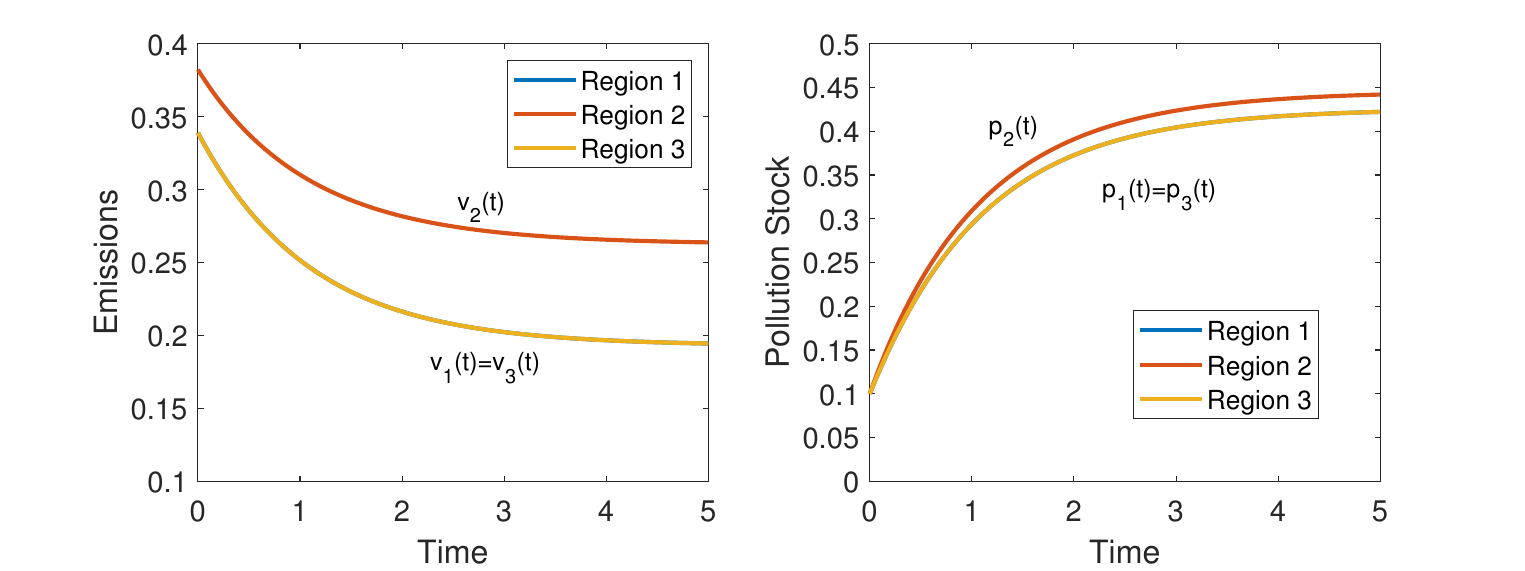}
\caption[Numerical emission (left) and pollution stock (right) time-paths along the equilibrium strategy obtained with the Chebyshev method.]{\label{exper2numsol} Numerical emission (left) and pollution stock (right) time-paths along the equilibrium strategy obtained with the Chebyshev method.}
\end{figure}

In Figure \ref{exper2numsol} we represent the emission (left) and pollution (right) time paths obtained with the Chebyshev numerical method. As expected, the optimal strategies and the pollution stocks of Players 1 and 3 coincide.

Unfortunately, for 3 or more players we lack of an explicit solution. Nevertheless,  we point out that, for the same values of $h, \text{TOL}$ and $N^p_i$, Chebyshev method outperforms the spline method in computational time cost.

\begin{figure}[h]
\centering
\includegraphics[width=7cm,height=5 cm]{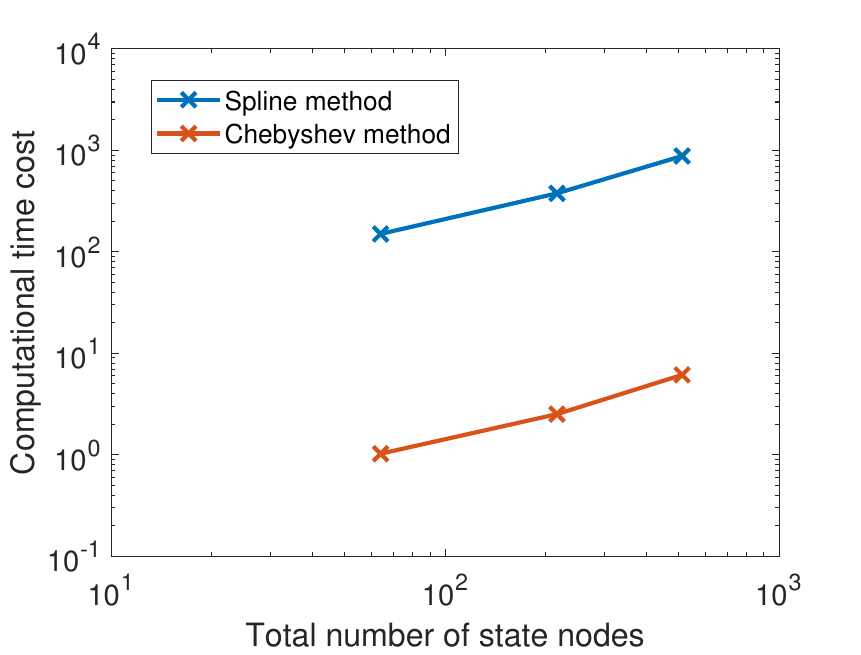}
\caption[Computational time cost of the spline(blue) and Chebyshev(red) methods for $N^p=3,5,7, \ h=10^{-3}, \ TOL=10^{-4}$.]{\label{exper2perfor} Computational time cost of the spline(blue) and Chebyshev(red) methods for $N^p=3,5,7, \ h=10^{-3}, \ TOL=10^{-4}$.}
\end{figure}

In Figure \ref{exper2perfor} we represent for the spline(blue) and Chebyshev(red) methods, the total number of spatial nodes $(N^p_i+1)^3$ vs the computational time cost for $N^p=3,5,7, \ h=10^{-3}, \ \text{TOL}=10^{-4}$. Other values for $h$ and $\text{TOL}$ were also tested, and the chosen ones are the fastest for the spline method.

For the same parameter values, the Chebyshev method requires, in average, 1/146 of the time of the spline method in order to obtain a numerical solution. This is not a complete performance analysis, since we lack of the explicit solutions, and we can not measure the numerical error. But point out that the results in the experiment for 2 players, and the fact that the qualitative solutions obtained with both methods are very similar, strongly suggest that the Chebyshev method outperforms the spline method.

\subsection{4 Players}

We now repeat Example 4 in \cite{deFrutos1}. The parameter values remain the same as in the previous experiment and the spatial configuration is given by

\begin{equation*}
K=[k_ij]=\left[\begin{matrix} -1 & 1 & 0 & 0 \\ 1 & -3 & 1 & 1 \\ 0 & 1 & -2 & 1 \\ 0 & 1 & 1 & -2 \end{matrix}\right]
\end{equation*}

This configuration means that Player 1 shares a frontier with Player 2, Player 2 shares a frontier with Players 1, 3, 4 and Player 3 shares a boundary with players 2 and 4. All the countries are isolated from outside. Under this configuration, Players 3 and 4 are ``symmetric'' since they share the same amount of frontiers with other countries and, therefore, their strategies should coincide.

\begin{figure}[h]
\centering
\includegraphics[width=13cm,height=5 cm]{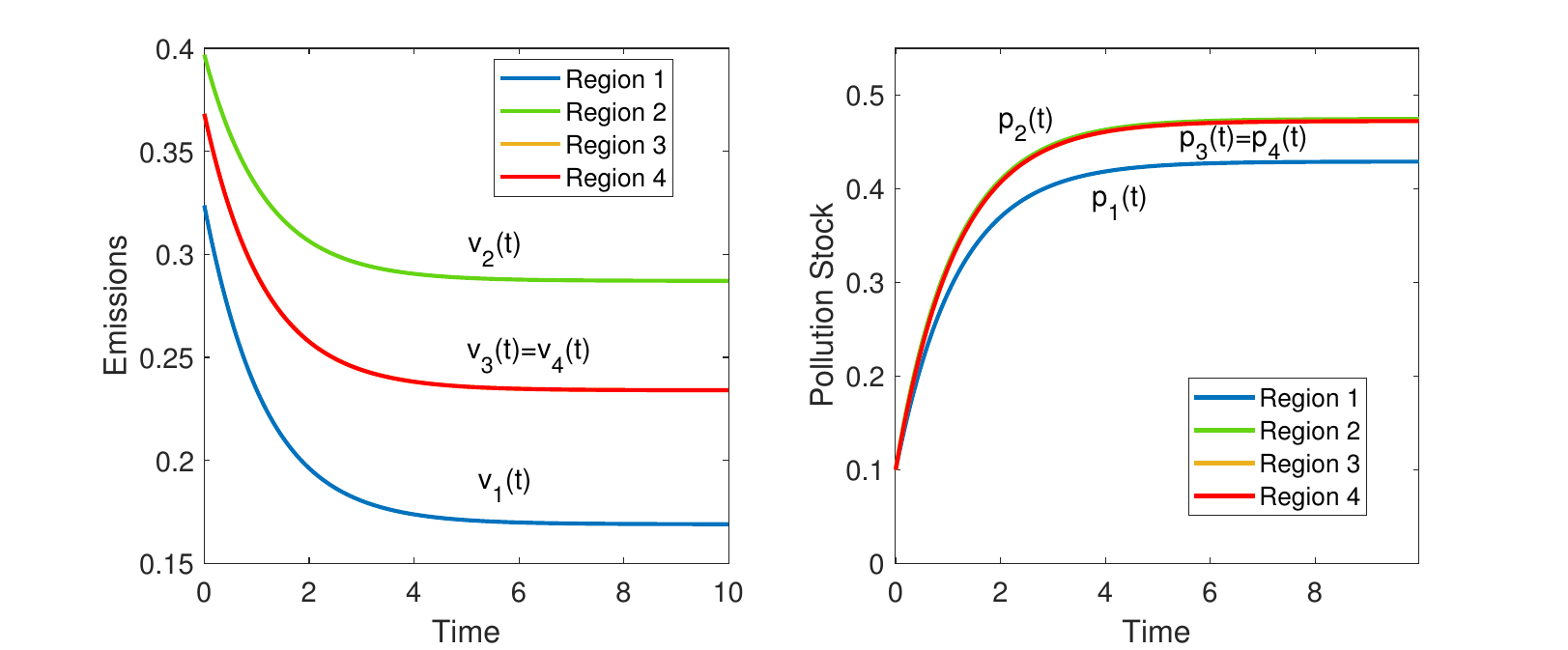}
\caption[Numerical emission (left) and pollution stock (right) time-paths along the equilibrium strategy obtained with the Chebyshev method.]{\label{exper3numsol} Numerical emission (left) and pollution stock (right) time-paths along the equilibrium strategy obtained with the Chebyshev method.}
\end{figure}

In Figure \ref{exper3numsol} we represent the emission (left) and pollution (right) time paths obtained with the Chebyshev numerical method. As expected, the optimal strategies and the pollution stock of Players 3 and 4 coincide.

Concerning numerical performance, the results are similar to the result in the experiment for 3 players. For the same values of $h, TOL$ and $N^p_i$, Chebyshev method outperforms in computational time cost the spline method.

\begin{figure}[h]
\centering
\includegraphics[width=7cm,height=5 cm]{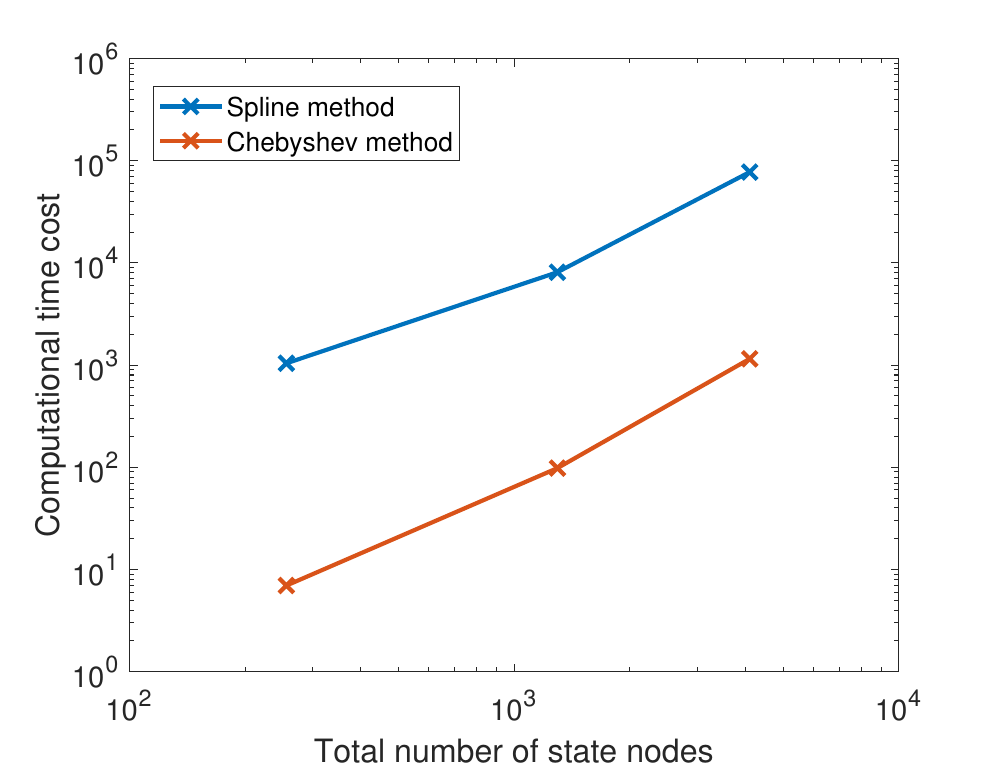}
\caption[Computational time cost of the spline(blue) and Chebyshev(red) methods for $N^p=3,5,7, \ dt=10^{-3}, \ TOL=10^{-4}$.]{\label{exper3perfor} Computational time cost of the spline(blue) and Chebyshev(red) methods for $N^p=3,5,7, \ dt=10^{-3}, \ TOL=10^{-4}$.}
\end{figure}

In Figure \ref{exper3perfor} we represent, for the spline(blue) and Chebyshev(red) methods, the total number of spatial nodes $(N^p_i+1)^4$ vs the computational time cost for $N^p=3,5,7, \ h=10^{-3}, \ TOL=10^{-4}$.

The Chebyshev method requires, in average, 1/100 of the time of the spline method in order to obtain a similar numerical solution.

As before, in the parallelization procedure, the optimal computational time cost is attained for a value $N_b$ such that $1<N_b<8^4=N_{\boldsymbol{P}}$.

\

Finally, we would like to point out that other experiments in \cite{deFrutos1}, including different spatial specifications and/or that one of the regions is not isolated from outside, have also been carried out. For not overloading the paper we have not included the results, but they have been similar to the ones presented in this work.

\section{Conclusions}\label{Conclus}

We have presented a tensorial-parallel Chebyshev collocation method for a game theory problem, which has a fairly good computational cost behaviour. This is due to the fact that it combines parallezation with some algorithms that allow, employing tensorization, to evaluate multidimensional Chebyshev polynomials efficently.

We should mention that the localization indexes presented (see Subsection \ref{multevalmulpol}) are not unique. Other dimension orders could be considered.

In this paper, we have presented the main ideas of a Chebyshev based algorithm which can be adapted to other differential game problems. These techniques may help to improve the numerical computation of problems which are affected by the known ``curse of dimensionality'', which appears when collocation methods are applied to problems with multiple dimensions.

Future work will be oriented in two different paths.

On one hand, in \cite{deFrutos2}, a Chebyshev based reduced function basis interpolation method is also presented. That technique allows to obtain the same numerical error with much less computational effort that a direct interpolation, as the one that we have employed in this work. Since the ``curse of dimensionality'' is still present, for a bigger number of players and number of state nodes, it would be interesting to adapt the reduced basis method to this problem.

On the other hand, we would like to adapt and test the algorithm to more complex model specifications. For example, it could be considered that each region $i$ can be divided in $n$ subregions, where player $i$ controls the emissions in each of the different subregions. Incorporate wind and a nonlinear reaction term in the pollution dynamics is also interesting since, although it is a model more computationally challenging, it is also closer to reality.

\subsection{Funding}

This research was supported by Junta de Castilla y Le\'{o}n cofinanced by FSE-YEI (first author) and Junta de Castilla y Le\'{o}n by project VA169P20 cofinanced by FEDER funds (second author).

\subsection{Acknowledgments}

The authors thank Javier de Frutos and Guiomar Mart\'{i}n-Herr\'{a}n for stimulating discussion.

\end{document}